\documentclass[12pt, reqno]{amsart}
\allowdisplaybreaks
\begin{document}
\setcounter{page}{1}

\title[\hfilneg \hfil Integrodifferential equations]
{Properties of Certain Partial dynamic Integrodifferential equations}

\author[D. B. Pachpatte\hfil \hfilneg]
{Deepak B. Pachpatte}

\address{Deepak B. Pachpatte \newline
 Department of Mathematics,
 Dr. Babasaheb Ambedkar Marathwada University, Aurangabad,
 Maharashtra 431004, India}
\email{pachpatte@gmail.com}


\begin{abstract}
 The aim of the present  paper is to study the existence, uniqueness and some other properties of solutions of a certain partial dynamic integrodifferential equations.The Banach fixed point theorem and certain fundamental inequality with explicit estimates are used to establish our results.
\end{abstract}

\maketitle

\section{Introduction}
 The study of time scale calculus was initiated by Stefan Hilger in his Ph.D dissertation which unifies the continuous and discrete calculus\cite{HIG}. Since then many authors have worked on various aspects dynamic equations on timescale calculus\cite{Dbp1,Dbp2,Dbp3,Dbp4,Dbp5}. Basic information on time scale calculus can be found in \cite{Boh1,Boh2,Boh3,HIG}.  Many authors have studied various types of  partial dynamic equations on time scales\cite{Dbp3,Dbp4,Jack,Sun,Men}.In \cite{Liu,Nowak,Xing} have studied the integrodifferential equations and its properties.  Motivated by the results in the above papers in this paper we study properties of certain partial dynamic integrodifferential equations.  
 In what follows $\mathbb{R}$ denotes the set of real numbers and $\mathbb{T}$ denotes the arbitrary time scales. Now we give some basic definitions of time scale calculus. The function $f:\mathbb{T} \to \mathbb{R}$ is said to be rd-continuous if $f$ is continuous at each right dense point of $\mathbb{T}$ and is denoted by $C_{rd}$. Let two time scales with at least two point be denoted by $\mathbb{T}_1$ and $\mathbb{T}_2$ and  $\Omega  = \mathbb{T}_1  \times \mathbb{T}_2 $. The delta partial derivative of a real valued function $f$ on $\mathbb{T}_1  \times \mathbb{T}_2 $ has a $\Delta _1 $ partial derivative $f^{\Delta _1 } \left( {t_1 ,t_2 } \right)$  with respect to $t_1$ if for each $\epsilon  > 0$ there exists a neighbourhood $U_{t_1 } $ of $t_1$ such that.
 \[
\left| {f\left( {\sigma _1 \left( {t_1 } \right),t_2 } \right) - f\left( {s,t_2 } \right) - f^{\Delta _1 } \left( {t_1 ,t_2 } \right)\left( {\sigma _1 \left( {t_1 } \right) - s} \right)} \right| \le \varepsilon \left| {\sigma _1 \left( {t_1 } \right) - s} \right|
\]
for all $s \in U_{t_2 }$.
The delta partial derivative of a real valued function $f$ on $\mathbb{T}_1  \times \mathbb{T}_2 $ has a $\Delta _2 $ partial derivative $f^{\Delta _1 } \left( {t_1 ,t_2 } \right)$  with respect to $t_2$ if for each $\eta  > 0$ there exists a neighbourhood $U_{t_2 } $ of $t_2$ such that
\[
\left| {f\left( {t_1 ,\sigma _2 \left( {t_2 } \right)} \right) - f\left( {t_1 ,l} \right) - f^{\Delta _2 } \left( {t_1 ,t_2 } \right)\left( {\sigma _2 \left( {t_2 } \right) - l} \right)} \right| \le \eta \left| {\sigma _2 \left( {t_2 } \right) - l} \right|
\]
 for all $u \in U_{t_1 }$.

The partial derivative of $w(x,y)$ for $\left( {x,y} \right) \in \Omega $ with respect to $x,y$ and $xy$ is denoted by $w^{\Delta _1 } \left( {x,y} \right), w^{\Delta _2 } \left( {x,y} \right)$ and $w^{\Delta _1 \Delta _2 } \left( {x,y} \right) = w^{\Delta _2 \Delta _1 } \left( {x,y} \right)$.
 Suppose $I=[a,b]$ with $a<b$ and $\overline \Omega =\Omega \times I$. The partial derivative of $u(x,y,z)$ for $(x,y,z) \in C_{rd} (\overline \Omega,R)$ with respect to $x,y$ and $xy$ is defined by $w^{\Delta _1 } \left( {x,y,z} \right),w^{\Delta _2 } \left( {x,y,z} \right)$ and $w^{\Delta _1 \Delta _2 } \left( {x,y,z} \right) = w^{\Delta _2 \Delta _1 } \left( {x,y,z} \right)$

 In this paper we study the partial dynamic integrodifferential equation of the form
 \[
u^{\Delta _2 \Delta _1 } \left( {x,y,z} \right) = F\left( {x,y,z,u(x,y,z),u^{\Delta _1 } (x,y,z),u^{\Delta _2 } (x,y,z),\left( {Hu} \right)(x,y,z)} \right),
 \tag{1.1}\]

with the conditions
 \[
u\left( {x,y_0 ,z} \right) = \alpha \left( {x,z} \right),\,u\left( {x_0 ,y,z} \right) = \beta \left( {y,z} \right)\,
\tag{1.2}\]
for $\left( {x,y} \right) \in \Omega $ where
\[
\left( {Hu} \right)\left( {x,y,z} \right) = \int\limits_a^b {G\left( {x,y,z,q,u\left( {x,y,q} \right),u^{\Delta _1 } \left( {x,y,q} \right),u^{\Delta _2 } \left( {x,y,q} \right)} \right)} \Delta q,
\tag{1.3}\]
where $G \in C_{rd} \left( {\overline \Omega \times \mathbb{R}^3 ,\mathbb{R}} \right),F \in C_{rd} \left( {\overline \Omega \times \mathbb{R}^4 ,\mathbb{R}} \right)$ and $\alpha ,\beta \,\, \in C_{rd} \left( {\mathbb{R}_ +   \times I,\mathbb{R}} \right)$.

We have
$u\left( {x_0 ,y_0 ,z} \right) = \alpha \left( {x_0 ,z} \right) = \beta \left( {x_0 ,z} \right)$.

Now for $u,u^{\Delta _1 } ,u^{\Delta _2 }  \in C_{rd} \left( {\overline \Omega,\mathbb{R}} \right)$, we denote
\[
\left| {u\left( {x,y,z} \right)} \right|_W  = \left| {u\left( {x,y,z} \right)} \right| + \left| {u^{\Delta _1 } \left( {x,y,z} \right)} \right| + \left| {u^{\Delta _2 } \left( {x,y,z} \right)} \right|.
\tag{1.4}\].

Let $S$ be the space function satisfying the condition
\[
\left| {u\left( {x,y,z} \right)} \right|_W  = O\left( {e_\lambda  \left( {x,y,\left| z \right|} \right)} \right),
\tag{1.5}\]
where $\lambda >0$ is a positive constant. In space $S$ we define norm $u$ by
\[
\left| u \right|_s  = \mathop {\sup }\limits_{\left( {x,y,z} \right) \in \Omega  \times I} \left[ {\left| {u\left( {x,y,z} \right)} \right|_w e_{\Theta \lambda } \left( {x,y,\left| z \right|} \right)} \right].
\tag{1.6}\]
The norm defined  $(1.6)$ is clearly a Banach Space.

Then $(1.5)$ implies that there is a constant $N \ge 0$ such that
\[
\left| {u\left( {x,y,z} \right)} \right|_w  \le N\left( {e_\lambda  \left( {x,y,\left| z \right|} \right)} \right),
\tag{1.7}\]
and we have
\[
\left| u \right|_s  \le N.
\tag{1.8}\]

The solution of $(1.1)$ and $(1.2)$ is a function  $u\left( {x,y,z} \right) \in C_{rd} \left( {\overline \Omega  ,\mathbb{R}^n } \right)$ satisfying  $(1.1)$ and $(1.2)$. It is easy to see that $u(x,y,z)$ with $(1.1)$ and $(1.2)$ satisfy the following dynamic integrodifferential equation.
\begin{align*}
&u\left( {x,y,z} \right) \\
 &= \alpha \left( {x,z} \right) + \beta \left( {y,z} \right) - \alpha \left( {0,z} \right) \\
 & + \int\limits_{x_0 }^x {\int\limits_{y_0 }^y {F\left( {s,t,z,u(s,t,z),u^{\Delta _1 } (s,t,z),u^{\Delta _2 } (s,t,z),\left( {Hu} \right)(s,t,z)} \right)\Delta t\Delta s} },
\tag{1.9}\end{align*}
for $(x,y,z) \in C_{rd}({\overline \Omega ,\mathbb{R}})$
\begin{align*}
 &u^{\Delta _1 } \left( {x,y,z} \right) \\
 &= \alpha ^{\Delta _1 } \left( {x,z} \right) \\
  &+ \int\limits_{y_0 }^y {F\left( {x,t,z,u(x,t,z),u^{\Delta _1 } (x,t,z),u^{\Delta _2 } (x,t,z),\left( {Hu} \right)(x,t,z)} \right)\Delta t},
\tag{1.10}
\end{align*}
\begin{align*}
 &u^{\Delta _2 } \left( {x,y,z} \right) \\
  &= \beta ^{\Delta _2 } \left( {y,z} \right) \\
   &+ \int\limits_{x_0 }^x {F\left( {s,y,z,u(x,t,z),u^{\Delta _1 } (s,y,z),u^{\Delta _2 } (s,y,z),\left( {Hu} \right)(s,y,z)} \right)\Delta s}.
\tag{1.11}
\end{align*}

We need following Lemma given in \cite{Boh3}.

Lemma [\cite{Boh3}, Theorem 2.6] Let $u \in C_{rd} \left( {\mathbb{T},\mathbb{R}_ +  } \right)$,$a \in \mathbb{R}_ + $
\[
u^\Delta  \left( t \right) \le a(t)u(t),
\]
for all $t \in \mathbb{T}^k $, then
\[
u(t) \le u(t_0 )e_a \left( {t,t_0 } \right),
\]
for all $t \in \mathbb{T}^k $.

\section{Main Results}
Now we give our main results

 \textbf{Theorem 1.1}  Suppose that the functions $F,G$ in $(1.1)$ satisfy the condition
 \begin{align*}
 &\left| {F\left( {x,y,z,u_1 ,u_2 ,u_3 ,u_4 } \right) - F\left( {x,y,z,\overline {u_1 } ,\overline {u_2 } ,\overline {u_3 } ,\overline {u_4 } } \right)} \right| \\
  &\le M\left( {x,y,z} \right)\left[ {\left| {u_1  - \overline {u_1 } } \right| + \left| {u_2  - \overline {u_2 } } \right| + \left| {u_3  - \overline {u_3 } } \right| + \left| {u_4  - \overline {u_4 } } \right|} \right],
 \tag{2.1}
 \end{align*}
 \begin{align*}
 & \left| {G\left( {x,y,z,q,u_1 ,u_2 ,u_3 } \right) - G\left( {x,y,z,q,\overline {u_1 } ,\overline {u_2 } ,\overline {u_3 } } \right)} \right| \\
 & \le K\left( {x,y,z,q} \right)\left[ {\left| {u_1  - \overline {u_1 } } \right| + \left| {u_2  - \overline {u_2 } } \right| + \left| {u_3  - \overline {u_3 } } \right|} \right],
   \tag{2.2}
  \end{align*}
 where $M \in C_{rd} \left( {\overline \Omega  ,\mathbb{R}_ +  } \right)$ and $K \in C_{rd} \left( {\overline \Omega  \times I ,\mathbb{R}_ +  } \right)$.

 For $\lambda$ as in $(1.5)$, there exists a nonnegative $\gamma _i \left( {i = 1,2,3} \right)$ such that
   \begin{align*}
 &\int\limits_{x_0 }^x {\int\limits_{y_0 }^y {M\left( {s,t,z} \right)} \left[ {e_\lambda  \left( {s,t,\left| z \right|} \right)} \right.}  \\
 &\left. { + \int\limits_a^b {k\left( {s,t,z,q} \right)e_\lambda  \left( {s,t,\left| q \right|} \right)\Delta q} } \right]\Delta t\Delta s \le \gamma _1 e_\lambda  \left( {x,y,\left| z \right|} \right),
  \tag{2.3}
 \end{align*}
\begin{align*}
 &\int\limits_{y_0 }^y {M\left( {x,t,z} \right)} \left[ {e_\lambda  \left( {x,t,\left| q \right|} \right)} \right. \\
  &\left. { + \int\limits_a^b {k\left( {x,t,z,q} \right)e_\lambda  \left( {x,t,\left| q \right|} \right)\Delta q} } \right]\Delta t \le \gamma _2 e_\lambda  \left( {x,y,\left| z \right|} \right),
   \tag{2.4}
\end{align*}
  \begin{align*}
 &\int\limits_{x_0 }^x {M\left( {s,y,z} \right)} \left[ {e_\lambda  \left( {x,t,\left| z \right|} \right)} \right. \\
 & \left. { + \int\limits_a^b {k\left( {s,y,z,q} \right)e_\lambda  \left( {s,y,\left| q \right|} \right)\Delta q} } \right]\Delta s \le \gamma _3 e_\lambda  \left( {x,y,\left| z \right|} \right),
  \tag{2.5}
\end{align*}
for $x,y \in \Omega$, $z \in I$.

 There exist nonnegative constants $\eta _i \left( {i = 1,2,3} \right)$ such that
 \begin{align*}
  &\left| {\alpha \left( {x,z} \right)} \right| + \left| {\beta \left( {y,z} \right)} \right| + \left| {\alpha \left( {0,z} \right)} \right| \\
  &+ \int\limits_{x_0 }^x {\int\limits_{y_0 }^y {\left| {f\left( {s,t,z,0,0,\left( {H0} \right)\left( {s,t,z} \right)} \right)} \right|} } \Delta t\Delta s \le \eta _1 e_\lambda  \left( {x,y,\left| z \right|} \right),
  \tag{2.6}
 \end{align*}
 \[
\left| {\alpha ^\Delta  \left( {x,z} \right)} \right| + \int\limits_{y_0 }^y {\left| {F\left( {x,t,z,0,0,0,\left( {H0} \right)\left( {x,t,z} \right)} \right)} \right|} \Delta t \le \eta _2 e_\lambda  \left( {x,y,\left| z \right|} \right),
\tag{2.7} \]
\[
\left| {\beta ^\Delta  \left( {y,z} \right)} \right| + \int\limits_{y_0 }^y {\left| {F\left( {s,y,z,0,0,0,\left( {H0} \right)\left( {s,y,z} \right)} \right)} \right|} \Delta s \le \eta _3 e_\lambda  \left( {x,y,\left| z \right|} \right),
\tag{2.8} \]
where $\alpha,\beta $  are as in $(1.2)$.

If $\gamma  = \gamma _1  + \gamma _2  + \gamma _3  < 1$ then problem $(1.1)-(1.2)$ has a unique solution $u(x,y,z)$ on $(1.1)-(1.2)$ in $S.$

 \textbf{Proof.} Let $u(x,y,z) \in S$ and define the operator $\mathbb{T}$ by
  \begin{align*}
 &\left( {Tu} \right)\left( {x,y,z} \right) \\
 &= \alpha \left( {x,z} \right) + \beta \left( {y,z} \right) - \alpha \left( {0,z} \right) \\
  &+ \int\limits_{x_0 }^x {\int\limits_{y_0 }^y {F\left( {s,t,z,u(s,t,z),u^{\Delta _1 } \left( {s,t,z} \right),u^{\Delta _2 } \left( {s,t,z} \right),\left( {Hu} \right)\left( {s,t,z} \right)} \right)} } \Delta t\Delta s.
  \tag{2.9}
 \end{align*}
Now we show that $P$ maps $S$ into itself. $Tu$ is rd-continuous on $\Omega \times I$ and $Tu \in R$.

From $(2.9)$ and given hypotheses we have
\begin{align*}
&\left( {Pu} \right)\left( {x,y,z} \right) \\
& \le \left| {\alpha \left( {x,z} \right)} \right| + \left| {\beta \left( {y,z} \right)} \right| + \left| {\alpha \left( {0,z} \right)} \right| \\
& + \int\limits_{x_0 }^x {\int\limits_{y_0 }^y {\left| {F\left( {s,t,z,u(s,t,z),u^{\Delta _1 } \left( {s,t,z} \right),u^{\Delta _2 } \left( {s,t,z} \right),\left( {Hu} \right)\left( {s,t,z} \right)} \right)} \right.} }  \\
&\left. { - F\left( {s,t,z,0,0,0,\left( {H0} \right)\left( {s,t,z} \right)} \right)} \right|\Delta t\Delta s \\
&+ \int\limits_{x_0 }^x {\int\limits_{y_0 }^y {\left| {F\left( {s,t,z,0,0,0,\left( {H0} \right)\left( {s,t,z} \right)} \right)} \right|} } \Delta t\Delta s \\
&  \le \eta _1 e_\lambda  \left( {x,y,\left| z \right|} \right) + \int\limits_{x_0 }^x {\int\limits_{y_0 }^y {M\left( {s,t,z} \right)} } \left[ {e_\lambda  \left( {s,t,\left| z \right|} \right)\left| {u\left( {s,t,z} \right)} \right|e_{\Theta \lambda } \left( {s,t,\left| z \right|} \right)} \right. \\
&\left. { + \int\limits_a^b {k\left( {x,y,z,q} \right)} e_\lambda  \left( {s,t,\left| q \right|} \right)\left| {u\left( {x,y,\left| q \right|} \right)} \right|_W e_{\Theta \lambda } \left( {s,t,\left| q \right|} \right)\Delta q} \right]\Delta t\Delta s \\
&\le \eta _1 e_\lambda  \left( {x,y,\left| z \right|} \right) + \left| u \right|_s \int\limits_{x_0 }^x {\int\limits_{y_0 }^y {M\left( {s,t,z} \right)} } \left[ {e_\lambda  } \right.\left( {s,t,\left| z \right|} \right) \\
&\left. { + \int\limits_a^b {k\left( {x,y,z,q} \right)} e_\lambda  \left( {s,t,\left| q \right|} \right)\Delta q} \right]\Delta t\Delta s \\
&\le \left[ {\eta _1  + N\gamma _1 } \right]e_\lambda  \left( {x,y,\left| z \right|} \right).
\tag{2.10}
\end{align*}
Delta differentiating on both sides of $(2.9)$ with respect to $x$ and  $(1.8)$ we have
\begin{align*}
&\left| {\left( {Pu} \right)^{\Delta _1 } \left( {x,y,z} \right)} \right| \\
&\le \alpha ^{\Delta _1 } \left( {x,z} \right) \\
&+ \int\limits_{y_0 }^y {\left| {F\left( {x,t,z,u(x,t,z),u^{\Delta _1 } \left( {x,t,z} \right),u^{\Delta _2 } \left( {x,t,z} \right),\left( {Hu} \right)\left( {x,t,z} \right)} \right)} \right.}  \\
&\left. { - F\left( {x,t,z,0,0,0,\left( {H0} \right)\left( {x,t,z} \right)} \right)} \right|\Delta t \\
&+ \int\limits_{y_0 }^y {\left| {F\left( {x,t,z,0,0,0,\left( {H0} \right)\left( {x,t,z} \right)} \right)} \right|} \Delta t \\
& \le \eta _2 e_\lambda  \left( {x,y,\left| z \right|} \right) + \left| u \right|_s \int\limits_{y_0 }^y {M\left( {x,t,z} \right)\left[ {e_\lambda  } \right.\left( {x,t,\left| z \right|} \right)}  \\
&\left. { + \int\limits_a^b {k\left( {x,t,z,q} \right)} e_\lambda  \left( {x,t,\left| q \right|} \right)\Delta q} \right]\Delta t \\
&\le \left[ {\eta _2  + N\gamma _2 } \right]e_\lambda  \left( {x,y,\left| z \right|} \right).
\tag{2.11}
\end{align*}
Similarly we have
\[
\left| {\left( {Pu} \right)^{\Delta _2 } \left( {x,y,z} \right)} \right| \le \left[ {\eta _3  + N\gamma _3 } \right]e_\lambda  \left( {x,y,\left| z \right|} \right).
\tag{2.12}\]
From $(2.10)-(2.12)$ we have
\[
\left| {Pu} \right|_s  \le \left[ {\left( {\eta _1  + \eta _2  + \eta _3 } \right) + N\gamma } \right].
\]
Thus proving that $P$ maps   $S$ into itself.

Now we show that operator $P$ is a contraction map. Let $u(x,y,z), \overline u \left( {x,y,z} \right) \in S $. From $(2.9)$ we have
\begin{align*}
&\left| {\left( {Pu} \right)\left( {x,y,z} \right) - \left( {P\overline u } \right)\left( {x,y,z} \right)} \right| \\
&\le \int\limits_{x_0 }^x {\int\limits_{y_0 }^y {\left| {F\left( {s,t,z,u(s,t,z),u^{\Delta _1 } (s,t,z),u^{\Delta _2 } (s,t,z),\left( {Hu} \right)\left( {s,t,z} \right)} \right)} \right.} }  \\
&\left. { - F\left( {s,t,z,\overline u (s,t,z),\overline u ^{\Delta _1 } (s,t,z),\overline u ^{\Delta _2 } (s,t,z),\left( {H\overline u } \right)\left( {s,t,z} \right)} \right)} \right|\Delta t\Delta s \\
&\le \left| {u - \overline u } \right|_s \int\limits_{x_0 }^x {\int\limits_{y_0 }^y {M\left( {x,t,z} \right)} } \left[ {e_\lambda  \left( {s,t,\left| z \right|} \right)} \right. \\ 
& \left. { + \int\limits_a^b {k\left( {s,t,z,q} \right)e_\lambda  \left( {s,t,\left| q \right|} \right)\Delta q} } \right]\Delta t\Delta s \\ 
&\le \left| {u - \overline u } \right|_s \gamma _1 e_\lambda  (x,y,\left| z \right|).
\tag{2.13}
\end{align*}

Similarly delta differentiating both sides of  $(2.12)$ with respect to  $x$ and $y$ we have
\[
\left| {\left( {Pu} \right)^{\Delta _1 } \left( {x,y,z} \right) - \left( {P\overline u } \right)^{\Delta _1 } \left( {x,y,z} \right)} \right| \le \left| {u - \overline u } \right|_s \gamma _2 e_\lambda  (x,y,\left| z \right|),
\tag{2.14} \]
and
\[
\left| {\left( {Pu} \right)^{\Delta _2 } \left( {x,y,z} \right) - \left( {P\overline u } \right)^{\Delta _2 } \left( {x,y,z} \right)} \right| \le \left| {u - \overline u } \right|_s \gamma _3 e_\lambda  (x,y,\left| z \right|).
\tag{2.15} \]
From $(2.13)-(2.15)$ we obtain
\[
\left| {Pu - P\overline u } \right|_s  \le \gamma \left| {u - \overline u } \right|_s.
\]
Since $\gamma <1 $, $P$ has a unique fixed point in $S$ by Banach fixed point theorem. The fixed point of $P$ is a solution of $(1.1)-(1.2)$. This completes the proof.

\section{Properties of solutions}
Now we study the properties of solution of dynamic integrodifferential equation of the form
\[
u^{\Delta _2 \Delta _1 } \left( {x,y,z} \right) = f\left( {x,y,z,u\left( {x,y,z} \right),\left( {hu} \right)\left( {x,y,z} \right)} \right),
\tag{3.1}\]
with $(1.2)$ for $(x,y,z) \in \overline \Omega $ where
\[
\left( {hu} \right)\left( {x,y,z} \right) = \int\limits_a^b {j\left( {x,y,z,q,u\left( {x,y,q} \right)} \right)} dq,
\tag{3.2}\]
in which $i \in C_{rd} \left( {\overline \Omega \times \mathbb{R},\mathbb{R}} \right)$, $f \in C_{rd} \left( {\Omega  \times \mathbb{R}^2 ,\mathbb{R}} \right)$.

Now we prove the following dynamic inequality which can be used in studying some properties of solutions.

 \textbf{Theorem 3.1}
Let $w,p \in C_{rd} \left( {\overline \Omega ,\mathbb{R}_ +  } \right)$, $\in C_{rd} \left( {\overline \Omega  \times I,\mathbb{R}_ +  } \right)$ and $c \ge 0$ a constant. If
\begin{align*}
 w\left( {x,y,z} \right)
 &\le c + \int\limits_{x_0 }^x {\int\limits_{y_0 }^y {\left[ {p\left( {s,t,z} \right)w\left( {s,t,z} \right)} \right.} }  \\
 &\left. { + \int\limits_a^b {r\left( {s,t,z,q} \right)w\left( {s,t,q} \right)\Delta q} } \right]\Delta t\Delta s,
\tag{3.3}\end{align*}
for $(x,y,z) \in \Omega $ then
\[
w\left( {x,y,z} \right) \le ce_{Q\left( {x,y,z} \right)} \left( {x,x_0 } \right),
\tag{3.4}\]
where $(x,y,z) \in \overline \Omega $ and 
\[
Q\left( {x,y,z} \right) = \int\limits_{y_0 }^y {\left[ {p\left( {s,t,Z} \right)\int\limits_a^b {r\left( {s,t,Z,q}  \right)\Delta q} } \right]}\Delta s.
\tag{3.5}\]

 \textbf{Proof.} For an arbitrary $Z \in I$ from $(3.3)$ we have

 \begin{align*}
 w\left( {x,y,Z} \right)
 &\le c + \int\limits_{x_0 }^x {\int\limits_{y_0 }^y {\left[ {p\left( {s,t,Z} \right)w\left( {s,t,Z} \right)} \right.} }\\
 &\left. { + \int\limits_a^b {r\left( {s,t,Z,q} \right)} w\left( {s,t,q} \right)\Delta q} \right]\Delta t\Delta s.
 \tag{3.6}
 \end{align*}
 Put
 \[
m\left( {s,t} \right) = p\left( {s,t,Z} \right)w\left( {s,t,Z} \right) + \int\limits_a^b {r\left( {s,t,z,q} \right)} w\left( {s,t,q} \right)\Delta q.
\tag{3.7}\]
The inequality $(3.6)$ becomes
 \[
w\left( {x,y,Z} \right) \le c + \int\limits_{x_0 }^x {\int\limits_{y_0 }^y {m(s,t)} } \Delta t\Delta s.
\tag{3.8}\]
Now define
\[
v\left( {x,y} \right) = c + \int\limits_{x_0 }^x {\int\limits_{y_0 }^y {m(s,t)} } \Delta t\Delta s,
\tag{3.9}\]
then
\[
v\left( {0,y} \right) = v\left( {x,0} \right) = c,w\left( {x,y,Z} \right) \le v\left( {x,y} \right).
\tag{3.10}\]
Delta differentiating both sides of $(3.9)$ with respect to $x$ and $y$ using $(3.7)$ and $(3.10)$ we have
\begin{align*}
 v^{\Delta _2 \Delta _1 } \left( {x,y} \right)
 & = m\left( {x,y} \right) \\
 & = p\left( {x,y,Z} \right)w\left( {x,y,Z} \right) + \int\limits_a^b {r\left( {x,y,Z,q} \right)} w\left( {x,y,q} \right)\Delta q \\
 & \le v\left( {x,y} \right)\left[ {p\left( {x,y,Z} \right) + \int\limits_a^b {r\left( {s,t,Z,q} \right)} \Delta q} \right].
\tag{3.11}
\end{align*}
By keeping $x$ fixed in $(3.11)$, and taking $y=t$ and delta integrating with respect to second variable from $y_0$ to $y$. Using the fact that
\begin{align*}
v^{\Delta _1 } \left( {x,y} \right)
&\le \int\limits_{y_0 }^y {\left[ {p\left( {x,t,Z} \right) + \int\limits_a^b {r\left( {x,t,Z,q} \right)} \Delta q} \right]} v\left( {x,t} \right)\Delta t \\
&\le v\left( {x,y} \right)\int\limits_{y_0 }^y {\left[ {p\left( {x,t,Z} \right) + \int\limits_a^b {r\left( {x,t,Z,q} \right)} \Delta q} \right]} \Delta t \\
&\le v\left( {x,y} \right)Q\left( {x,y,Z} \right).
\tag{3.12}
\end{align*}
 Now treating $y$ fixed in $(3.12)$ and applying Lemma we have
 \[
v(x,y) \le ce_{Q\left( {x,y,Z} \right)} \left( {x,x_0 } \right).
\tag{3.13}
\]
Because $Z$ is arbitrary and using $(3.10)$ we get $(3.9)$.

  \textbf{Theorem 3.2} Suppose the functions $f,j$ in  $(3.1)$,$(3.2)$ satisfy the conditions
  \[
\left| {f\left( {x,y,z,u,v} \right) - f\left( {x,y,z,\overline u ,\overline v } \right)} \right| \le p_1 \left( {x,y,z} \right)\left[ {\left| {u - \overline u } \right| + \left| {v - \overline v } \right|} \right],
\tag{3.14}\]
\[
\left| {j\left( {x,y,z,q,u} \right) - j\left( {x,y,z,\overline q ,\overline u } \right)} \right| \le p_2 \left( {x,y,z,q} \right)\left| {u - \overline u } \right|,
\tag{3.15}\]
where $p_1  \in C_{rd} \left( {\overline \Omega ,\mathbb{R}_ +  } \right)$, $p_2  \in C_{rd} \left( {\overline \Omega  \times I,\mathbb{R}_ +  } \right)$, $c \ge 0$
and
\[
\int\limits_{x_0 }^x {\int\limits_{y_0 }^y {\left[ {p_1 \left( {s,t,z} \right) + \int\limits_a^b {p_1 \left( {s,t,z,q} \right)\Delta q} } \right]} } \Delta t\Delta s < \infty,
\tag{3.16}\]
then the problem $(3.1)-(1.1)$ has at most one solution.

\textbf{Proof.}
  Let $u_1(x,y,z)$ and $u_2(x,y,z)$ be two solutions of problem $(3.1)-(1.1)$.
\begin{align*}	
 &\left| {u_1 \left( {x,y,z} \right) - u_2 \left( {x,y,z} \right)} \right| \\
&\le \int\limits_{x_0 }^x {\int\limits_{y_0 }^y {\left| {f\left( {s,t,z,u_1 \left( {s,t,z} \right),\left( {hu_1 } \right)\left( {s,t,z} \right)} \right)} \right.} }  \\
&\left. { - f\left( {s,t,z,u_2 \left( {s,t,z} \right),\left( {hu_2 } \right)\left( {s,t,z} \right)} \right)} \right|\Delta t\Delta s \\
& \le \int\limits_{x_0 }^x {\int\limits_{y_0 }^y {\left[ {p_1 \left( {s,t,z} \right)\left| {u_1 \left( {s,t,z} \right) - u_2 \left( {s,t,z} \right)} \right|} \right.} }  \\
&\left. { + \left| {\left( {hu_1 } \right)\left( {s,t,z} \right) - \left( {hu_2 } \right)\left( {s,t,z} \right)} \right|} \right]\Delta t\Delta s \\
& \le \int\limits_{x_0 }^x {\int\limits_{y_0 }^y {\left[ {p_1 \left( {s,t,z} \right)\left| {u_1 \left( {s,t,z} \right) - u_2 \left( {s,t,z} \right)} \right|} \right.} }  \\
&\left. { + \int\limits_a^b {p_1 \left( {s,t,z,q} \right)\left| {u_1 \left( {s,t,q} \right) - u_2 \left( {s,t,q} \right)} \right|\Delta q} } \right]\Delta t\Delta s.
\tag{3.17}
\end{align*}
Now applying Theorem $3.1$to $(3.17)$  yields $\left| {u_1 \left( {x,y,z} \right) - u_2 \left( {x,y,z} \right)} \right| \le 0$ which gives
$u_1 \left( {x,y,z} \right) = u_2 \left( {x,y,z} \right)$. This proves that there is at most one solution to problem $(3.1)-(1.1)$.

Now we prove the theorem which gives the boundedness of solution of $(3.1)-(1.1)$.

 \textbf{Theorem 3.3.} Suppose the function $f,j,\alpha, \beta$ in $(3.1)-(1.1)$ satisfy the conditions
\[
\left| {f\left( {x,y,z,u,v} \right)} \right| \le p_1 \left( {x,y,z} \right)\left[ {\left| u \right| + \left| v \right|} \right],
\tag{3.18}\]
\[
\left| {j\left( {x,y,z,u,v} \right)} \right| \le p_2 \left( {x,y,z,q} \right)\left| u \right|,
\tag{3.19}\]
\[
\left| {\alpha \left( {x,z} \right) + \beta \left( {y,z} \right) - \alpha \left( {0,z} \right)} \right| \le c,
\tag{3.20}\]
where $p_1  \in C_{rd} \left( {\Omega ,\mathbb{R}_ +  } \right)$, $p_2  \in C_{rd} \left( {\Omega  \times I,\mathbb{R}_ +  } \right)$, $c \ge 0$ is a constant and the condition $(3.16)$ holds. Then solution $u(x,y,z)$ is bounded and
\[
\left| {u\left( {x,y,z} \right)} \right| \le ce_{Q\left( {x,y,z} \right)} \left( {x,x_0 } \right),
\tag{3.21}\]
for $(x,y,z) \in \overline \Omega $

\textbf{Proof.} Since $u(x,y,z)$ is a solution of $(3.1)-(1.1)$. We have
\begin{align*}
\left| {u\left( {x,y,z} \right)} \right|
&  \le \left| {\alpha \left( {x,z} \right) + \beta \left( {y,z} \right) - \alpha \left( {0,z} \right)} \right| \\
& + \int\limits_{x_0 }^x {\int\limits_{y_0 }^y {\left| {f\left( {s,t,z,u_{} \left( {s,t,z} \right),\left( {hu} \right)\left( {s,t,z} \right)} \right)} \right|} } \Delta t\Delta s \\
&\le c + \int\limits_{x_0 }^x {\int\limits_{y_0 }^y {\left[ {p_1 \left( {s,t,z} \right)\left| {u\left( {s,t,z} \right)} \right|} \right.} }  \\
&\left. { + \int\limits_a^b {p_2 \left( {s,t,z,q} \right)\left| {u\left( {s,t,q} \right)} \right|\Delta q} } \right]\Delta t\Delta s.
\tag{3.22}
\end{align*}

Now an application of Theorem $3.1$ to $(3.22)$ yields $(3.21)$ thus proving the boundedness of solution.

Now we give the dependency of solution of equation on given condition

\textbf{Theorem 3.4.} Suppose the function $f,k$ in $(3.1),(3.2)$ satisfy the conditions $(3.14),(3.15)$ and the condition $(3.16)$ holds. Let $u(x,y,z)$ and $v(x,y,z)$ be the solutions of equation with condition $(1.2)$ and
\[
v\left( {x,0,z} \right) = \overline \alpha  \left( {x,z} \right),\,\,v\left( {0,y,z} \right) = \overline \beta  \left( {y,z} \right),
\tag{3.23}\]
respectively and
\[
\left| {\alpha \left( {x,z} \right) + \beta \left( {y,z} \right) - \alpha \left( {0,z} \right) - \left[ {\overline \alpha  \left( {x,z} \right) + \overline \beta  \left( {y,z} \right) - \overline \alpha  \left( {0,z} \right)} \right]} \right| \le a,
\tag{3.24}\]
where $\alpha,\beta, \overline \alpha, \overline \beta \in C_{rd}(\mathbb{R}_+ \times I, \mathbb{R}) $ and $a \ge 0$ is constant. Then
\[
| u( x,y,z )  -v(x,y,z)| \le ae_{Q\left( {x,y,z} \right)} \left( {x,x_0 } \right).
\tag{3.25}\]

\textbf{Proof.} Since $u(x,y,z)$ and $v(x,y,z)$ are solutions of $(3.1)$-$(1.1)$ and $(3.1)-(3.23)$ and the given conditions we have
 \begin{align*}
&\left| {u\left( {x,y,z} \right) - u\left( {x,y,z} \right)} \right| \\
&\le \left| {\alpha \left( {x,z} \right) + \beta \left( {y,z} \right) - \alpha \left( {0,z} \right)} \right. \\
&\left. { - \left[ {\overline \alpha  \left( {x,z} \right) + \overline \beta  \left( {y,z} \right) - \overline \alpha  \left( {0,z} \right)} \right]} \right| \\
&+ \int\limits_{x_0 }^x {\int\limits_{y_0 }^y {\left| {f\left( {s,t,z,u\left( {s,t,z} \right),\left( {hu} \right)\left( {s,t,z} \right)} \right)} \right.} }  \\
&\left. { - f\left( {s,t,z,v\left( {s,t,z} \right),\left( {hu} \right)\left( {s,t,z} \right)} \right)} \right|\Delta t\Delta s \\
&\le a + \int\limits_{x_0 }^x {\int\limits_{y_0 }^y {\left[ {p_1 \left( {s,t,z} \right)\left| {u\left( {s,t,z} \right) - v\left( {s,t,z} \right)} \right|} \right.} }  \\
&
 \left. { + \int\limits_a^b {p_2 \left( {s,t,z,q} \right)\left| {u\left( {s,t,z} \right) - v\left( {s,t,z} \right)} \right|\Delta q} } \right]\Delta t\Delta s.
\tag{3.26}
\end{align*}

Now an application of Theorem $3.1$ to $(3.26)$ gives the estimate $(3.25)$ which gives the dependency of solution of equation $(3.1)$ on given conditions.

 \textbf{Acknowledgement.} This research is supported by Science and Engineering Research Board(SERB), New Delhi, India, Sanct. No. SB/S4/MS:861/13

\end{document}